\def\E{\end{document}}
\documentclass[11pt]{article}
\usepackage{amssymb,amsmath}

\begin{document}
\title{\bf
A  quasilinear Schr\"{o}dinger equation with Hartree type nonlinearity}
  \author{Xianfa Song$^{a}${\thanks{Corresponding author, E-mail: songxianfa2004@163.com(X.F. Song)
 }} \ and
 Zhi-Qiang Wang$^{b,c}${\thanks{Email:zhi-qiang.wang@usu.edu(Z. Q. Wang)}}\\
\small a  Department of Mathematics, School of Mathematics, Tianjin University,
\small Tianjin, 300072, P. R. China\\
\small b College of Mathematics and Informatics, Fujian Normal University,
\small Fuzhou 350117, P. R. China\\
\small c  Department of Mathematics and Statistics, Utah State University,
\small  Logan, UT 84322, USA
}

\maketitle
\date{}

\newtheorem{theorem}{Theorem}[section]
\newtheorem{definition}{Definition}[section]
\newtheorem{lemma}{Lemma}[section]
\newtheorem{proposition}{Proposition}[section]
\newtheorem{corollary}{Corollary}[section]
\newtheorem{remark}{Remark}[section]
\renewcommand{\theequation}{\thesection.\arabic{equation}}
\catcode`@=11 \@addtoreset{equation}{section} \catcode`@=12

\begin{abstract}

In this paper, we deal with the Cauchy problem of the quasilinear Sch\"{o}dinger equation
\begin{equation*}
\left\{
\begin{array}{lll}
iu_t=\Delta u+2uh'(|u|^2)\Delta h(|u|^2)+(W(x)\ast|u|^2)u,\ x\in \mathbb{R}^N,\ t>0\\
u(x,0)=u_0(x),\quad x\in \mathbb{R}^N.
\end{array}\right.
\end{equation*}
Here $h(s)$ and $W(x)$ are some real valued functions. Our focus is to investigate how the interplay between the potential $W(x)$ and the quasilinear presence $h(s)$ affects the blowup in finite time and global existence of the solution. In a special, we can obtain the watershed condition on $W(x)$ in the following sense: If $W(x)\in L^1(\mathbb{R}^N)\cap \{L^q(\mathbb{R}^N)+L^{\infty}(\mathbb{R}^N)\} $, then exist $q_c$ and $q_s$ such that the solution is global existence for any initial data in the energy space when $q>q_c$ and the solution maybe
 blow up in finite time for some initial data when $q_s<q<q_c$, and for $q=q_c$ whether the solution is global existence or not depend on the initial data.



\end{abstract}

\section{Introduction}
\qquad In this paper, we consider the following Cauchy problem:
\begin{equation}
\label{1} \left\{
\begin{array}{lll}
iu_t=\Delta u+2uh'(|u|^2)\Delta h(|u|^2)+(W(x)\ast|u|^2)u, \ x\in \mathbb{R}^N, \ t>0\\
u(x,0)=u_0(x),\quad x\in \mathbb{R}^N.
\end{array}\right.
\end{equation}
Here $N\geq 3$, $h(s)$ and $W(x)$ are some real functions, $h(s)\geq 0$ for $s\geq 0$, and $W(x)$ is even. We always assume $h$ is a smooth function. Quasilinear equation  (\ref{1}) with various choices of $h$ can be used to model a lot of physical phenomena, such as the self-channelling of a high-power ultra short laser in matter. It often appears in plasma physics and dissipative quantum mechanics, and in condensed matter theory, see \cite{Bass, Borovskii,  Bouard, Goldman, Litvak, Makhankov, Ritchie}.
By the classic results in Section 3.2 of \cite{Cazenave}, if the even function $W(x)\in L^1(\mathbb{R}^N)\cap \{L^q(\mathbb{R}^N)+L^{\infty}(\mathbb{R}^N)\}$ for some $q>1$, then $(W*|u|^2)u\in C(H^1(\mathbb{R}^N), H^{-1}(\mathbb{R}^N))$, $(W*|u|^2)u\in C(L^r(\mathbb{R}^N), L^{r'}(\mathbb{R}^N))$ and $(W(x)\ast |u|^2)u\in C_b(\mathbb{R}^N, B_M(0))$, where $r=\frac{4q}{2q-1}$ and $B_M(0)=\{u\in \mathbb{C}:|u|\leq M\}$. Consequently, the local well-posedness result on (\ref{1}) can be contained in these frame work of \cite{Colin1, Kenig, Poppenberg1}.

In this paper, we focus on the effect of the potential $W(x)$ and we investigate quantitatively how the interplay between the potential and the quasilinear term influent the blowup in finite time and global existence of the solution of (\ref{1}).
The definition of the global existence and blowup in finite time for the solution of (\ref{1}) is given below.

{\bf Definition 1.} {\it Let $u(x,t)$ be the solution of (\ref{1}). We say that $u(x,t)$ exists globally if the maximal existence interval for $t$ is $[0, +\infty)$, and we say $u(x,t)$ blows up in finite time if there exists a finite time $T>0$ such that
\begin{align}
\lim_{t\rightarrow T^-} \int_{\mathbb{R}^N}[|\nabla u(x,t)|^2+|\nabla h(|u(x,t)|^2)|^2)]dx=+\infty.
\end{align}
 }

The motivations of our study in this paper are as follows.

First of all,  the Cauchy problem of semilinaer Schr\"{o}dinger equations with Hartree type nonlinearity has been studied extensively in the literature, which dealt with the global existence and other behaviors of the solutions. We can refer to \cite{Cao11, Cho13, Cho17, Ivanov08, Zagatti92} and the references therein. The study for quasilinear equations like would be natural extentions of these work, and here we study the properties for the solution of quasilinear equation (\ref{1}) which also contains Hartree type nonlinearity in the equation. This is an area that has not been studied in the past.

Secondly, there are many works in the literature about the global existence and blowup phenomena of semilinear Schr\"{o}dinger equation. The following Cauchy problem
\begin{equation}
\label{2} \left\{
\begin{array}{lll}
iu_t=\Delta u+F(|u|^2)u \quad {\rm for} \ x\in \mathbb{R}^N, \ t>0\\
u(x,0)=u_0(x),\quad x\in \mathbb{R}^N
\end{array}\right.
\end{equation}
 was considered by Glassey in the celebrated  work \cite{Glassey}. The key condition on the blowup of the solution to (\ref{2}) is  that $sF(s)\geq c_N G(s)$ for some constant $c_N>1+\frac{2}{N}$ and all $s\geq 0$, here $G(s)=\int_0^s F(t)dt$. In \cite{Berestycki}, Berestycki and Cazenave established a sharp threshold on the blowup of the solution. Other related results can be found in \cite{Cazenave, Zhang1} and the references therein. In \cite{Guo}, the Cauchy problem of quasilinear Sch\"{o}dinger equation
\begin{equation}
\label{xj} \left\{
\begin{array}{lll}
-iu_t+\Delta u+2(\Delta |u|^2)u+|u|^{q-2}u=0\quad {\rm for} \ x\in \mathbb{R}^N, \ t>0\\
u(x,0)=u_0(x),\quad x\in \mathbb{R}^N
\end{array}\right.
\end{equation}
was studied by Guo, Chen and Su. They showed that the solution of (\ref{xj}) blows up in finite time if $4+\frac{4}{N}<q<2\cdot 2^*$ for some initial data.
 Letting
\begin{align}
F(x,|u|^2)=[W(x)\ast |u|^2],\label{08911}
\end{align}
(\ref{1}) can be written in the following form
\begin{equation}
\label{1'} \left\{
\begin{array}{lll}
iu_t=\Delta u+2uh'(|u|^2)\Delta h(|u|^2)+F(x,|u|^2)u \quad {\rm for} \ x\in \mathbb{R}^N, \ t>0\\
u(x,0)=u_0(x),\quad x\in \mathbb{R}^N.
\end{array}\right.
\end{equation}
In \cite{SW1}, the authors of this paper gave qualitative analysis for (\ref{1'}) with $F(x,s)=F(s)$ and obtained the key conditions on the global existence and blowup in finite time for the solution, which are the explicit relationships between $sF(s)$ and $G(s)$(for example, $c_NG(s)\leq sF(s)$), the explicit relationships between $G(s)$(or $sF(s)$) and $h(s)$(for example, $[G(s)]^{q_i}\leq c [h(s)]^{2^*}$), where $G(s)=\int_0^s F(\eta)d\eta$, $2^*=\frac{2N}{N-2}$. However, letting $\tilde{G}(x,s)=\int_0^s F(x,\eta)d\eta$, one cannot find the explicit relationship between $sF(x,s)$ and $\tilde{G}(x,s)$.
We need to develop a new approach for establishing the key condition on the global existence and blowup in finite time for the solution.

The last motivation is inspired by the results on the following problem:
\begin{equation}
\label{mhs2} \left\{
\begin{array}{lll}
iu_t=\Delta u+2b^2\alpha |u|^{2\alpha-2}u \Delta (|u|^{2\alpha})+(p(x)*|u|^2)u \quad {\rm for} \ x\in \mathbb{R}^N\setminus \{\mathbf{0}\}, \ t>0\\
u(x,0)=u_0(x),\quad x\in \mathbb{R}^N.
\end{array}\right.
\end{equation}
Here
\begin{equation}
\label{55w1} p(x)=p(r)=\left\{
\begin{array}{lll}
\frac{1}{r^p},\quad 0<r\leq 1\\
f(r),\quad 1\leq r\leq 2\\
\frac{1}{r^C},\quad r\geq 2, \ r=|x|
\end{array}\right.
\end{equation}
$C>\max(N,p)$, and $f(r)>0$ satisfies $pf(r)+rf'(r)\leq 0\leq Cf(r)+rf'(r)$ if $1\leq r\leq 2$.
Very recently, we found that: The solution of (\ref{mhs2}) is global existence for any $u_0$ which belongs some energy space if $p<p_c$ while the solution will blow up in finite time under certain conditions if $p\geq p_c$. Here $p_c=\frac{N[\max(2\alpha,1)\cdot 2^*-2]}{2^*}$ if $0<\alpha<\frac{N-1}{N}$. These results are also parts of this paper. Naturally, we wish to generalize the results on (\ref{mhs2}) to more general cases.  In fact, it is the first time that we let the potential $W(x)$ be the criterion of the conditions on the blowup in finite time and global existence of the solution to a quasilinear Sch\"{o}dinger equation.

We use $C_s$ to denote the best constant in the Sobolev's inequality
\begin{align}
\int_{\mathbb{R}^N}|w|^{2^*}dx\leq C_s\left(\int_{\mathbb{R}^N}|\nabla w|^2dx\right)^{\frac{2^*}{2}}\quad {\rm for \ any}\quad w\in H^1(\mathbb{R}^N),\label{zjcs}
\end{align}

Our first result will establish some sufficient conditions on the global existence of the solution to (\ref{1}).

{\bf Theorem 1.} {\it  Let $u(x,t)$ be the solution of (\ref{1}) with $u_0\in X$,
\begin{align}
X=\{w\in H^1(\mathbb{R}^N),\quad \int_{\mathbb{R}^N}|\nabla h(|w|^2)|^2dx<+\infty\}.\label{kongjianziji}
\end{align}
Assume that there exist $a>0$ in the sense of infimum and $\alpha>0$ in the sense of supremum such that $\max(s^{\frac{1}{2}}, s^{\alpha})\leq a[h(s)+s^{\frac{1}{2}}]$ for $s\geq 1$.\\
Assume $W(x)$ is even and satisfies
$$
(C1)\ W(x)=W_1(x)+W_2(x)\in L^1(\mathbb{R}^N)\cap \{ L^q(\mathbb{R}^N)+L^{\infty}(\mathbb{R}^N)\},\ q>1,\ q>\frac{\max(\alpha,\frac{1}{2})2^*}{\max(2\alpha, 1) 2^*-2}.
$$
Then we have the following assertions.

(1). If $W(x)\leq 0$ for all $x\in \mathbb{R}^N$, then $u$ is global existence for any $u_0\in X$.

(2). Assume $W(x)\geq 0$ for all $x\in\mathbb{R}^N$ or changes sign. We have three subcases.

Case (i) If $0<\alpha<\frac{N-1}{N}$, if $q>\frac{2^*}{\max(2\alpha, 1) 2^*-2}$, then the solution is global existence for any initial $u_0\in X$;

Case (ii) If $0<\alpha<\frac{N-1}{N}$, if $q=\frac{2^*}{\max(2\alpha, 1)2^*-2}$, then the solution is global existence for the initial $u_0\in X$ satisfying
$$
a^22^{\frac{(q-1)N+2q}{qN}}C_s^{\frac{2}{2^*}}\|W\|_{L^q(\mathbb{R}^N)}\|u_0\|_{L^2(\mathbb{R}^N)}^{\min(4-4\alpha,2)}< 1;
$$

Case (iii) If $\alpha\geq\frac{N-1}{N}$, if $q>1$, then the solution is global existence for any initial $u_0\in X$.
}

{\bf Remark 1.1.} 1. If $h(|u|^2)\equiv 0$, (\ref{1}) becomes the classic semilinear Schr\"{o}dinger equation with Hartree type
nonlinearity. Then we can take $a=1$ and $\alpha=\frac{1}{2}$, the solution of (\ref{1}) is global existence if $q>\frac{N}{2}$.


2. If $h(|u|^2)=b|u|^{2\alpha}(b\geq 0)$, $0<\alpha<\frac{N-1}{N}$ and $W(x)=p(x)$, where $p(x)$ is defined as (\ref{55w1}), then the solution is global existence when $p<\frac{N[\max(2\alpha,1)\cdot 2^*-2]}{2^*}$.

3. We suspect that the assumption $W(x)\in L^1(\mathbb{R}^N)\cap \{L^q(\mathbb{R}^N)+L^{\infty}(\mathbb{R}^N)\}$ in this paper can be weaken as
$W(x)\in \{L^q(\mathbb{R}^N)+L^{\infty}(\mathbb{R}^N)\}$. However, by the local well-posedness result on (\ref{1}) in \cite{Poppenberg1}, we need $W(x)\in L^1(\mathbb{R}^N)$ to ensure that $(W(x)\ast |u|^2)u\in C_b(\mathbb{R}^N, B_M(0))$.

Our second result is about sufficient conditions on the blowup in finite time for the solution of (\ref{1}).

{\bf Theorem 2.} {\it Let $u(x,t)$ be the solution of (\ref{1}) with $u_0\in X$, $xu_0\in L^2(\mathbb{R}^N)$, $E(u_0)\leq 0$ and $\Im \int_{\mathbb{R}^N}\bar{u}_0(x\cdot \nabla u_0)dx>0$. Assume that there exists constant $k$ such that for $s\geq 0$, $sh''(s)\leq kh'(s)$ if $h'(s)\geq 0$ or  $sh''(s)\geq kh'(s)$ if $h'(s)\leq 0$. Assume further
$
W(x)\geq 0,\quad W(x)\in L^1(\mathbb{R}^N)\cap \{ L^q(\mathbb{R}^N)+L^{\infty}(\mathbb{R}^N)\}\quad {\rm for \ some }\quad q>1.
$
Suppose the following holds.
$$(C2)\qquad \qquad \quad \quad[\max((2k+1)N,0)+2]W+(x\cdot \nabla W)\leq 0.\qquad \qquad \qquad$$
Then there exists a finite time $T$ such that
$$
\lim_{t\rightarrow T^-} \int_{\mathbb{R}^N}[|\nabla u(x,t)|^2+|\nabla h(|u|^2)(x,t)|^2]dx=+\infty.
$$
}

{\bf Remark 1.2.} It is easy to check that for $s\geq 0$, $sh''(s)\leq kh'(s)$ if $h'(s)\geq 0$ (or $sh''(s)\geq kh'(s)$ if $h'(s)\leq 0$) implies that $-h'(s)h''(s)s\geq -k(h'(s))^2$. We would like to discuss
conditions $(C2)$ and $E(u_0)\leq 0$in details.

1. Obviously, $W(x)\equiv c>0$ means that $[\max((2k+1)N,0)+2]W+(x\cdot \nabla W)>0$, which is the opposite of the condition (C2).
$W(x)\equiv c<0$ implies that $E(u_0)>0$ for $u_0(x)\neq 0$.

2. It is well known that $L^{q_1}(\mathbb{R}^N)\cap L^{q_2}(\mathbb{R}^N)\neq \emptyset$ for $0<q_1<q_2$.
If $W(x)$ is a nontrivial radially symmetric function, by the condition $$W(r)\geq 0,\quad [\max((2k+1)N,0)+2]W(r)+rW'(r)\leq 0,$$ then for any $q\geq\frac{N}{\max[(2k+1)N,0]+2}$, $${\rm if}\quad 0<r<1,\quad  W(r)\geq \frac{W(1)}{r^{[\max((2k+1)N,0)+2]}}\notin L^1(\mathbb{R}^N)\cap \{L^q(\mathbb{R}^N)+L^{\infty}(\mathbb{R}^N)\}.$$ However, $\frac{1}{r^{[\max((2k+1)N,0)+2]}}\in L^{\tilde{q}}(\mathbb{R}^N)+L^{\infty}(\mathbb{R}^N)$ for  $\tilde{q}<\frac{N}{\max[(2k+1)N,0]+2}$.

To illustrates how the potential $W(x)$ affects the blowup in finite time and global existence of the solution, we give a corollary below.
Consider
\begin{equation}
\label{mhs3} \left\{
\begin{array}{lll}
iu_t=\Delta u+2\alpha |u|^{2\alpha-2}u \Delta (|u|^{2\alpha})+(W*|u|^2)u \quad {\rm for} \ x\in \mathbb{R}^N, \ t>0\\
u(x,0)=u_0(x),\quad x\in \mathbb{R}^N.
\end{array}\right.
\end{equation}
Here we assume $0<\alpha<\frac{N-1}{N}$.  Define
$$
q_c=\frac{2^*}{\max(2\alpha, 1)2^*-2}.
$$
{\bf Corollary 1.1.} {\it Assume $W(x)\in L^1(\mathbb{R}^N)\cap \{L^q(\mathbb{R}^N)+L^{\infty}(\mathbb{R}^N)\}$ and let $u$ be the solution of
(\ref{mhs3}). We have the following assertions.\\
(1) If  $q>q_c$, then the solution is global existence for any initial data $u_0\in X$.\\
(2)
 If $W(x)\in L^1(\mathbb{R}^N)\cap \{L^{q_c}(\mathbb{R}^N)+L^{\infty}(\mathbb{R}^N)\}$,  then the solution is global existence for $u_0$ satisfying $$2^{\frac{(q_c-1)N+2q_c}{q_cN}}C_s^{\frac{2}{2^*}}\|W\|_{L^{q_c}(\mathbb{R}^N)}\|u_0\|_{L^2(\mathbb{R}^N)}^{\min(4-4\alpha,2)}=1.$$
(3)  If $W(x)\geq 0$, $W(x)\in L^1(\mathbb{R}^N)\cap \{L^q(\mathbb{R}^N)+L^{\infty}(\mathbb{R}^N)\}$ for some $q_c>q>1$ and $[(2\alpha-1)N+2]W+(x\cdot \nabla W)\leq 0$, then the solution will blow up in finite time for $u_0\in X$, $xu_0\in L^2(\mathbb{R}^N)$, $E(u_0)\leq 0$ and $\Im \int_{\mathbb{R}^N}\bar{u}_0(x\cdot \nabla u_0)dx>0$.
}

{\bf Remark 1.3.} By Remark 1.2, the conditions $W(x)\geq 0$ and $[(2\alpha-1)N+2]W+(x\cdot \nabla W)\leq 0$ imply that
$W(x)\not\in  L^1(\mathbb{R}^N)\cap \{L^q(\mathbb{R}^N)+L^{\infty}(\mathbb{R}^N)\}$ for any $q>q_c$, yet $W(x)\in L^1(\mathbb{R}^N)\cap \{L^q(\mathbb{R}^N)+L^{\infty}(\mathbb{R}^N)\}$ may hold for some $1<q\leq q_c$. In another word, the function $W(x)\in L^1(\mathbb{R}^N)\cap \{L^q(\mathbb{R}^N)+L^{\infty}(\mathbb{R}^N)\}$ can be classified and be taken as criterion of the blowup and global existence of the solution of (\ref{1}) in the following sense: The solution of (\ref{mhs3}) is always global existence if $W(x)\in L^1(\mathbb{R}^N)\cap \{L^q(\mathbb{R}^N)+L^{\infty}(\mathbb{R}^N)\}$ for $q>q_c$, but the solution may blow up in finite time if $W(x)\in L^1(\mathbb{R}^N)\cap \{L^q(\mathbb{R}^N)+L^{\infty}(\mathbb{R}^N)\}$
 for some $1<q\leq q_c$.

The organization of this paper is as follows. In Section 2, we will prove the mass and energy conservation laws and some equalities.
In Section 3, we will prove Theorem 1.  In Section 4, we will prove Theorem 2. In Section 5, we will construct a sharp threshold for the global existence and blowup in finite time for the solution of (\ref{1}). In Section 6, we will establish the pseudoconformal conservation law and consider the asymptotic behaviour for the solution.

\section{Preliminaries}
\qquad In convenience, we will use $C$, $C'$, and so on, to denote some constants in the sequels, the values of it may vary line to line.
We define the mass and energy of (\ref{1}) as follows.

(i) Mass: $$ m(u)=\left(\int_{\mathbb{R}^N}|u(\cdot,t)|^2dx\right)^{\frac{1}{2}}:=[M(u)]^{\frac{1}{2}};$$

(ii) Energy : $$E(u)=\frac{1}{2}\int_{\mathbb{R}^N}[|\nabla u|^2+|\nabla h(|u|^2)|^2]dx-\frac{1}{4}\int_{\mathbb{R}^N}(W*|u|^2)|u|^2dx.$$

We will give a lemma about conservations of these quantities as follows.

{\bf Lemma 2.1.} {\it Assume that $u$ is the solution to (\ref{1}). Then in the time interval $[0,t]$ when it exists, $u$ satisfies

(i) Mass conversation: $$ m(u)=\left(\int_{\mathbb{R}^N}|u(x,t)|^2dx\right)^{\frac{1}{2}}=\left(\int_{\mathbb{R}^N}|u_0(x)|^2dx\right)^{\frac{1}{2}}=m(u_0);$$

(ii) Energy conversation: $$E(u)=\frac{1}{2}\int_{\mathbb{R}^N}[|\nabla u|^2+|\nabla h(|u|^2)|^2]dx-\frac{1}{4}\int_{\mathbb{R}^N}(W*|u|^2)|u|^2dx=E(u_0);$$

(iii) $$\frac{d}{dt} \int_{\mathbb{R}^N}|x|^2|u|^2dx=-4\Im \int_{\mathbb{R}^N} \bar{u}(x\cdot \nabla u)dx;$$

(iv) \begin{align}
&\quad \frac{d}{dt} \Im \int_{\mathbb{R}^N} \bar{u}(x\cdot \nabla u)dx=-2\int_{\mathbb{R}^N}|\nabla u|^2dx-(N+2)\int_{\mathbb{R}^N}|\nabla h(|u|^2)|^2dx\nonumber\\
&\qquad -8N\int_{\mathbb{R}^N}h''(|u|^2)h'(|u|^2)|u|^4|\nabla u|^2dx-\frac{1}{2}\int_{\mathbb{R}^N}[(x\cdot \nabla W)*|u|^2]|u|^2dx.\label{10131}
\end{align}

}

{\bf Proof:} (i) Multiplying (\ref{1}) by $2\bar{u}$, taking the imaginary part of the result, we get
\begin{align}
\frac{\partial }{\partial t}|u|^2=\Im(2\bar{u}\Delta u) =\nabla \cdot (2\Im \bar{u}\nabla u).\label{10121}
\end{align}
Integrating it over $\mathbb{R}^N\times [0,t]$, we have
$$ \int_{\mathbb{R}^N}|u|^2dx=\int_{\mathbb{R}^N}|u_0|^2dx.$$
Consequently, $m(u)=m(u_0)$.

(ii)  Multiplying (\ref{1}) by $2\bar{u}_t$, taking the real part of the result, then integrating it over $\mathbb{R}^N\times [0,t]$, we obtain
\begin{align*}
&\quad\int_{\mathbb{R}^N}[|\nabla u|^2+|\nabla h(|u|^2)|^2]dx-\frac{1}{2}\int_{\mathbb{R}^N}(W*|u|^2)|u|^2dx\\
&=\int_{\mathbb{R}^N}[|\nabla u_0|^2+|\nabla h(|u_0|^2)|^2]dx-\frac{1}{2}\int_{\mathbb{R}^N}(W*|u_0|^2)|u_0|^2dx.
\end{align*}
Therefore, $E(u)=E(u_0)$.

(iii) Multiplying (\ref{10121}) by $|x|^2$ and integrating it over $\mathbb{R}^N$, we get
\begin{align*}
\frac{d}{dt}\int_{\mathbb{R}^N}|x|^2|u|^2dx&=\int_{\mathbb{R}^N}|x|^2\nabla \cdot(2\Im (\bar{u}\nabla u))dx
=-4\Im \int_{\mathbb{R}^N}\bar{u}(x\cdot \nabla u)dx.
\end{align*}

(iv) Let  $a(x,t)=\Re u(x,t)$ and $b(x,t)=\Im u(x,t)$. Then
\begin{align*}
\int_{\mathbb{R}^N}\frac{d}{dt}\Im \bar{u}(x\cdot \nabla u)dx&=\int_{\mathbb{R}^N}\sum_{k=1}^N[x_k(b_t)_{x_k}a-x_k(a_t)_{x_k}b]dx+\int_{\mathbb{R}^N}\sum_{k=1}^N(x_kb_{x_k}a_t-x_ka_{x_k}b_t)dx\nonumber\\
&=\int_{\mathbb{R}^N} \sum_{k=1}^N[x_k(b_t)_{x_k}a-x_k(a_t)_{x_k}b]dx
+\int_{\mathbb{R}^N} \sum_{k=1}^N (x_ka_{x_k}\Delta a+x_kb_{x_k}\Delta b)dx\nonumber\\
&\quad+\frac{1}{2}\int_{\mathbb{R}^N}\sum_{k=1}^N x_k(|u|^2)_{x_k}[2h'(|u|^2)\Delta h(|u|^2)+(W*|u|^2)]dx\nonumber\\
&=N\int_{\mathbb{R}^N}(a_tb-ab_t)dx+\int_{\mathbb{R}^N}\sum_{k=1}^N(x_kb_{x_k}a_t-x_ka_{x_k}b_t)dx\nonumber\\
&\quad+\frac{N-2}{2}\int_{\mathbb{R}^N}|\nabla u|^2dx+\frac{N-2}{2}\int_{\mathbb{R}^N}|\nabla h(|u|^2)|^2dx\nonumber\\
&\quad-\frac{1}{2}\int_{\mathbb{R}^N}[(NW+\frac{(x\cdot \nabla W)}{2})*|u|^2]|u|^2dx,
\end{align*}
and
\begin{align*}
\frac{d}{dt}\int_{\mathbb{R}^N}\Im \bar{u}(x\cdot \nabla u)dx&=N\int_{\mathbb{R}^N}\left([a\Delta a+b\Delta b]+2|u|^2h'(|u|^2)\Delta h(|u|^2)+|u|^2(W*|u|^2)\right)dx\nonumber\\
&\quad +(N-2)\int_{\mathbb{R}^N}|\nabla u|^2dx+(N-2)\int_{\mathbb{R}^N}|\nabla h(|u|^2)|^2dx\nonumber\\
&\quad-\int_{\mathbb{R}^N}[(NW+\frac{(x\cdot \nabla W)}{2})*|u|^2]|u|^2dx\nonumber\\
&=-2\int_{\mathbb{R}^N}|\nabla u|^2dx-(N+2)\int_{\mathbb{R}^N}|\nabla h(|u|^2)|^2dx\nonumber\\
&-8N\int_{\mathbb{R}^N}h'(|u|^2)h''(|u|^2)|u|^4|\nabla u|^2dx-\int_{\mathbb{R}^N}[\frac{(x\cdot \nabla W)}{2}*|u|^2]|u|^2dx.
\end{align*}

{\bf Remark 2.1.} If $u$ is the solution of (\ref{1}), similar to \cite{SW1}, we have
\begin{align*}
\int_{\mathbb{R}^N}|u|^{p_2}dx&\leq \left(\int_{\mathbb{R}^N}|u|^2dx\right)^{\frac{p_1-p_2}{p_1-2}}\left(\int_{\mathbb{R}^N}|u|^{p_1}dx\right)^{\frac{p_2-2}{p_1-2}}\\
&=\left(\int_{\mathbb{R}^N}|u_0|^2dx\right)^{\frac{p_1-p_2}{p_1-2}}\left(\int_{\mathbb{R}^N}|u|^{p_1}dx\right)^{\frac{p_2-2}{p_1-2}}
\end{align*}
for $p_1>p_2>2$ by mass conservation law.

\section{The proof of Theorem 1}

\qquad In this section, we will prove Theorem 1 and establish the sufficient conditions on the global existence of the solution to (\ref{1}).

{\bf Proof of Theorem 1:} Under the assumptions of Theorem 1, we know that $\int_{\mathbb{R}^N}(|W|*|u_0|^2)|u_0|^2dx$ and $E(u_0)$ are well defined in $X$.

Case (1). $W(x)\leq 0$ for $x\in \mathbb{R}^N$. The global existence of the solution is a direct result of the energy conversation law of Lemma 2.1(ii) because
$$
\int_{\mathbb{R}^N}|\nabla u|^2dx+\int_{\mathbb{R}^N}|\nabla h(|u|^2)|^2dx\nonumber\\
+\frac{1}{2}\int_{\mathbb{R}^N}(|W|*|u|^2)|u|^2dx=2E(u_0),
$$
which implies that $\int_{\mathbb{R}^N}|\nabla u|^2dx+\int_{\mathbb{R}^N}|\nabla h(|u|^2)|^2dx+\int_{\mathbb{R}^N}(|W|*|u|^2)|u|^2dx$ is uniformly bounded for all $t>0$.

Case (2).  $W(x)\geq 0$ for all $x\in\mathbb{R}^N$ or changes sign.

Obviously, $q>\frac{\max(\alpha,\frac{1}{2})2^*}{\max(2\alpha, 1) 2^*-2}$ is equivalent to $\frac{4q}{2q-1}<\max(2\alpha, 1) 2^*$.

Denoting
\begin{align}
\tau_1=\frac{(2q-1)[\max(2\alpha,1)2^*-2]}{(2q-1)[\max(2\alpha,1)2^*-2]-2},\quad \tau_2=\frac{(2q-1)[\max(2\alpha,1)2^*-2]}{2},\label{117w1}
\end{align}
then
\begin{align}
\frac{2q-1}{q\tau_1}=\frac{(2q-1)[\max(2\alpha,1)2^*-2]-2}{q[\max(2\alpha,1)2^*-2]}, \quad \frac{2q-1}{q\tau_2}=\frac{2}{q[\max(2\alpha,1)2^*-2]}.\label{117w2}
\end{align}

Suppose that $W(x)=W_1(x)+W_2(x)\in L^1(\mathbb{R}^N)\cap \{L^q(\mathbb{R}^N)+L^{\infty}(\mathbb{R}^N)\}$, where $W_1(x)\in L^q(\mathbb{R}^N)$ and $W_2(x)\in L^{\infty}(\mathbb{R}^N)$.
Denote
$$
C_W=\frac{1}{2}\left(\int_{\mathbb{R}^N}|W_1|^qdx\right)^{\frac{1}{q}}.
$$

Using the mass and energy conversation laws of Lemma 2.1(ii), using H\"{o}der's inequality, Young's inequality, then Sobolev's inequality, we have
\begin{align}
&\quad\int_{\mathbb{R}^N}|\nabla u|^2dx+\int_{\mathbb{R}^N}|\nabla h(|u|^2)|^2dx=2E(u_0)+\frac{1}{2}\int_{\mathbb{R}^N}(W*|u|^2)|u|^2dx\nonumber\\
&=2E(u_0)+\frac{1}{2}\int_{\mathbb{R}^N}(W_1*|u|^2)|u|^2dx+\frac{1}{2}\int_{\mathbb{R}^N}(W_2*|u|^2)|u|^2dx\nonumber\\
&\leq C+\frac{1}{2}\left(\int_{\mathbb{R}^N}|W_1|^qdx\right)^{\frac{1}{q}}\left(\int_{\mathbb{R}^N}|u|^{\frac{4q}{2q-1}}dx\right)^{\frac{2q-1}{q}}
+\frac{1}{2}\|W_2\|_{L^{\infty}}\left(\int_{\mathbb{R}^N}|u|^2dx\right)^2\nonumber\\
&\leq C+\frac{1}{2}\|W_2\|_{L^{\infty}}(m(u_0))^4+2^{\frac{q-1}{q}}C_W\left(\int_{\{|u|\leq 1\}}|u|^{\frac{4q}{2q-1}}dx\right)^{\frac{2q-1}{q}}\nonumber\\
&\qquad+2^{\frac{q-1}{q}}C_W\left(\int_{\{|u|>1\}}|u|^{\frac{4q}{2q-1}}dx\right)^{\frac{2q-1}{q}}\nonumber\\
&\leq C+\frac{1}{2}\|W_2\|_{L^{\infty}}(m(u_0))^4+2^{\frac{q-1}{q}}C_W\left(\int_{\{|u|\leq 1\}}|u|^2dx\right)^{\frac{2q-1}{q}}\nonumber\\
&\qquad +2^{\frac{q-1}{q}}C_W\left(\int_{\{|u|>1\}}|u|^2dx\right)^{\frac{2q-1}{q\tau_1}}
\left(\int_{\{|u|>1\}}|u|^{\max(2\alpha,1)2^*}dx\right)^{\frac{2q-1}{q\tau_2}}\nonumber\\
&\leq C+\frac{1}{2}\|W_2\|_{L^{\infty}}(m(u_0))^4+2^{\frac{q-1}{q}}C_W\left(\int_{\mathbb{R}^N}|u|^2dx\right)^{\frac{2q-1}{q}}\nonumber\\
&\qquad+2^{\frac{q-1}{q}}C_W\left(\int_{\mathbb{R}^N}|u|^2dx\right)^{\frac{2q-1}{q\tau_1}}\left(\int_{\{|u|>1\}}|u|^{\max(2\alpha,1)2^*}dx\right)^{\frac{2q-1}{q\tau_2}}\nonumber\\
&\leq C'+2^{\frac{q-1}{q}}C_W\left(\int_{\mathbb{R}^N}|u_0|^2dx\right)^{\frac{2q-1}{q\tau_1}}\left(\int_{\{|u|>1\}}\max(|u|^{2\alpha\cdot2^*},|u|^{2^*})dx\right)^{\frac{2q-1}{q\tau_2}}\nonumber\\
&\leq C'+2^{\frac{q-1}{q}}C_W\left(\int_{\mathbb{R}^N}|u_0|^2dx\right)^{\frac{2q-1}{q\tau_1}}\left(\int_{\{|u|>1\}}a^{2^*}[h(|u|^2)+|u|]^{2^*}dx\right)^{\frac{2q-1}{q\tau_2}}\nonumber\\
&\leq C'+a^{\frac{2^*(2q-1)}{q\tau_2}}2^{\frac{(2^*-1)(2q-1)}{q\tau_2}}2^{\frac{q-1}{q}}C_W\|u_0\|_{L^2}^{\frac{2(2q-1)}{q\tau_1}}\left(\int_{\mathbb{R}^N}([h(|u|^2)]^{2^*}+|u|^{2^*})dx\right)^{\frac{2q-1}{q\tau_2}}\nonumber\\
&\leq C'+a^{\frac{2^*(2q-1)}{q\tau_2}}(C_s)^{\frac{(2q-1)}{q\tau_2}}2^{\frac{(2q-1)N+2q}{qN}}C_W\|u_0\|_{L^2}^{\frac{2(2q-1)}{q\tau_1}}\left(\int_{\mathbb{R}^N}[|\nabla h(|u|^2)|^2+|\nabla u|^2]dx\right)^{\frac{2^*(2q-1)}{2q\tau_2}}.\label{9141}
\end{align}

We give the estimate for (\ref{9141}) in three subcases.

Subcase (i) $0<\alpha< \frac{N-1}{N}=\frac{2^*+2}{22^*}$ and $q>\frac{2^*}{\max(2\alpha, 1) 2^*-2}$. In this subcase,
$$
\frac{2^*(2q-1)}{2q\tau_2}<1.
$$
 For any initial data $u_0$, using Young's inequality, we have
$$
\int_{\mathbb{R}^N}|\nabla u|^2dx+\int_{\mathbb{R}^N}|\nabla h(|u|^2)|^2dx\leq C+C'(C_W,u_0,q)+\frac{1}{2}\int_{\mathbb{R}^N}[|\nabla h(|u|^2)|^2+|\nabla u|^2]dx,
$$
which implies that
$$\int_{\mathbb{R}^N}|\nabla u|^2dx+\int_{\mathbb{R}^N}|\nabla h(|u|^2)|^2dx\leq C.$$

Subcase (ii) $0<\alpha<\frac{N-1}{N}=\frac{2^*+2}{22^*}$ and $q=\frac{2^*}{\max(2\alpha, 1) 2^*-2}$. In this subcase, $$\frac{(2q-1)}{q\tau_1}=2-\max(2\alpha,1)=\min(2-2\alpha,1),\qquad \frac{2^*(2q-1)}{2q\tau_2}=1.$$ For the initial $u_0$ satisfying
$$
a^{\frac{2^*(2q-1)}{q\tau_2}}(2^{2^*-1}C_s)^{\frac{(2q-1)}{q\tau_2}}2^{\frac{q-1}{q}}C_W\|u_0\|_{L^2}^{\frac{2(2q-1)}{q\tau_1}}<1,
$$
i.e.,
$$a^22^{\frac{(q-1)N+2q}{qN}}C_s^{\frac{2}{2^*}}\|W_1\|_{L^q}\|u_0\|_2^{\min(4-4\alpha,2)}<1,$$
we have
$$\left(1-a^22^{\frac{(q-1)N+2q}{qN}}C_s^{\frac{2}{2^*}}\|W_1\|_{L^q}\|u_0\|_2^{\min(4-4\alpha,2)}\right)\int_{\mathbb{R}^N}[|\nabla h(|u|^2)|^2+|\nabla u|^2]dx\leq C.$$

Subcase (iii) $\alpha\geq \frac{N-1}{N}=\frac{2^*+2}{22^*}$. In this subcase, $$\frac{2^*}{\max(2\alpha, 1) 2^*-2}\leq 1, \quad \frac{\max(\alpha,\frac{1}{2})2^*}{[\max(2\alpha,1)2^*-2]}<1.$$
Then for any $q>1$, we have
$$
\frac{4q}{2q-1}<\max(2\alpha, 1) 2^*,\quad \frac{2^*(2q-1)}{2q\tau_2}<1.
$$
Similar the proof of Subcase (i), applying Young's inequality to (\ref{9141}), we get $$\int_{\mathbb{R}^N}|\nabla u|^2dx+\int_{\mathbb{R}^N}|\nabla h(|u|^2)|^2dx\leq C$$
for all $t>0$.\hfill $\Box$

{\bf Remark 3.1.} 1.  In Subcase (i) and Subcase (ii)  above, if $0<\alpha\leq \frac{2}{2^*}$, $N\geq 4$, then $$\frac{2^*}{\max(2\alpha, 1) 2^*-2}>\frac{\max(\alpha,\frac{1}{2})2^*}{[\max(2\alpha,1)2^*-2]}\geq 1,$$ while if $\frac{2}{2^*}<\alpha<\frac{2^*+2}{22^*}=\frac{N-1}{N}$, $N\geq 4$, then
 $$\frac{2^*}{\max(2\alpha, 1) 2^*-2}>1>\frac{\max(\alpha,\frac{1}{2})2^*}{[\max(2\alpha,1)2^*-2]}.$$
If $0<\alpha<\frac{N-1}{N}$, $N=3$, then $$
\frac{2^*}{\max(2\alpha, 1) 2^*-2}>1>\frac{\max(\alpha,\frac{1}{2})2^*}{\max(2\alpha,1)2^*-2}.
$$

 2. The conclusion of (2)(iii) shows the interaction between the term $2uh'(|u|^2)\Delta h(|u|^2)$ and Hartree type nonlinearity. Roughly, if $h(s)$ increases fast enough when $s>1$, then the solution is global existence for any $W(x)\in L^1(\mathbb{R}^N)\cap \{L^q(\mathbb{R}^N)+L^{\infty}(\mathbb{R}^N)\}$($q>1$) and initial data $u_0\in X$.

\section{The proof of Theorem 2}
\qquad In this section, we will give the proof of Theorem 2 and deal with the sufficient conditions on blowup in finite time for the solution by using the results of Lemma 2.1.

{\bf Proof of Theorem 2:} Wherever $u$ exists, let
$$
y(t)=\Im \int_{\mathbb{R}^N}\bar{u}(x\cdot \nabla u)dx.
$$

We discuss it in two cases:

Case 1. $h(s)\equiv 0$ or $h(s)\neq 0$ and $k\leq -\frac{1}{2}$. We have

\begin{align}
\dot{y}(t)&=-2\int_{\mathbb{R}^N}|\nabla u|^2dx-(N+2)\int_{\mathbb{R}^N}|\nabla h(|u|^2)|^2dx-8N\int_{\mathbb{R}^N}h''(|u|^2)h'(|u|^2)|u|^4|\nabla u|^2dx\nonumber\\
&\quad -\int_{\mathbb{R}^N}[\frac{(x\cdot \nabla W)}{2}*|u|^2]|u|^2dx\nonumber\\
&\geq -2\int_{\mathbb{R}^N}|\nabla u|^2dx-(N+2+2kN)\int_{\mathbb{R}^N}|\nabla h(|u|^2)|^2dx-\int_{\mathbb{R}^N}[\frac{(x\cdot \nabla W)}{2}*|u|^2]|u|^2dx\nonumber\\
&=-4E(u_0)-(2k+1)N\int_{\mathbb{R}^N}|\nabla h(|u|^2)|^2dx-\frac{1}{2}\int_{\mathbb{R}^N}([2W+(x\cdot \nabla W)]*|u|^2)|u|^2dx\nonumber\\
&\geq 0,\label{10141}
\end{align}
which means that $y(t)\geq y(0)>0$ for $t>0$.

Case 2. $h(s)\neq 0$ and $k>-\frac{1}{2}$. We have
\begin{align}
\dot{y}(t)&\geq -2\int_{\mathbb{R}^N}|\nabla u|^2dx-(N+2+2kN)\int_{\mathbb{R}^N}|\nabla h(|u|^2)|^2dx-\int_{\mathbb{R}^N}[\frac{(x\cdot \nabla W)}{2}*|u|^2]|u|^2dx\nonumber\\
&=(2k+1)N\int_{\mathbb{R}^N}|\nabla u|^2dx-2[(2k+1)N+2]E(u_0)\nonumber\\
&\quad-\frac{1}{2}\int_{\mathbb{R}^N}([((2k+1)N+2)W+(x\cdot \nabla W)]*|u|^2)|u|^2dx\nonumber\\
&\geq 0,\label{10141}
\end{align}
which also means that $y(t)\geq y(0)>0$ for $t>0$.

Setting
$$
J(t)=\int_{\mathbb{R}^N}|x|^2|u|^2dx,
$$
we have $J'(t)=-4y(t)<-4y(0)<0$. Then
$$0\leq J(t)=J(0)+\int_0^tJ'(\tau)d\tau<J(0)-4y(0)t,$$
which implies that the maximum existence interval of time for $u$ is finite, and $u$ will blow up before $\frac{J(0)}{4y(0)}$.\hfill $\Box$

Now we give the proof of Corollary 1.1.

{\bf The proof of Corollary 1.1:} The results on the global existence of the solution are the direct results of Theorem 1. Since $k=\alpha-1$ and
$$\max((2k+1)N,0)+2=\max((2\alpha-1)N,0)+2,$$
the assumptions of Theorem 2 are satisfied, the solution will blow up in finite time for initial data $u_0$ satisfying $E(u_0)<0$ and $\Im \int_{\mathbb{R}^N}\bar{u}(x\cdot \nabla u_0)dx>0$.

We give another Corollary below, which establishes the conditions on blowup in finite time and global existence of the solution to (\ref{mhs2}).

{\bf Corollary 4.1.} {\it Assume that $u$ be the solution of (\ref{mhs2}). We have

Case 1. $b=0$.  The solution is global existence for any initial data $u_0\in H^1(\mathbb{R}^N)$ if $0<p<2$, while the solution will blow up in finite time for initial data $u_0$ satisfying $E(u_0)<0$ and $\Im \int_{\mathbb{R}^N}\bar{u}(x\cdot \nabla u_0)dx>0$ if $p\geq 2$.

Case 2. $b\neq 0$.  The solution is global existence for any $u_0\in X$ in one of the following subcases: Subcase (i) $0<\alpha<\frac{N-1}{N}$ and $p<\frac{N[\max(2\alpha,1)2^*-2]}{2^*}$; Subcase (ii) $\alpha\geq \frac{N-1}{N}$ and $p<N$. While the solution will blow up in finite time for $u_0$ satisfying $E(u_0)<0$ and $\Im \int_{\mathbb{R}^N}\bar{u}(x\cdot \nabla u_0)dx>0$ if $\alpha>0$ and $p\geq\frac{N[\max(2\alpha,1)2^*-2]}{2^*}$.
}

{\bf Proof:} By the definition of  $W(x)$, the condition $$[\max((2k+1)N,0)+2]W+(x\cdot \nabla W)\leq 0$$ implies that $p\geq \max((2k+1)N,0)+2$.

Case 1. If $b\equiv 0$, we take $k=-\frac{1}{2}$ and get $p\geq 2$. The solution will blow up in finite time for initial data $u_0$ satisfying $E(u_0)<0$ and $\Im \int_{\mathbb{R}^N}\bar{u}(x\cdot \nabla u_0)dx>0$ if $p\geq 2$. While $p<2$, then $W(x)\in L^1(\mathbb{R}^N)\cap \{L^q(\mathbb{R}^N)+L^{\infty}(\mathbb{R}^N)\}$ for some $q>\frac{N}{2}$, the solution is global existence for any initial data $u_0\in H^1(\mathbb{R}^N)$.

Case 2. If $b\neq 0$, we can take $k=\alpha-1$ and get
$$p\geq \max((2k+1)N,0)+2=\max((2\alpha-1)N,0)+2=\frac{N[\max(2\alpha,1)2^*-2]}{2^*}.$$
The solution will blow up in finite time for initial data $u_0$ satisfying $E(u_0)<0$ and $\Im \int_{\mathbb{R}^N}\bar{u}(x\cdot \nabla u_0)dx>0$ if $\alpha>0$ and $p\geq\frac{N[\max(2\alpha,1)2^*-2]}{2^*}$ by the results of Theorem 2. Subcase (i) $0<\alpha<\frac{N-1}{N}$ and $p<\frac{N[\max(2\alpha,1)2^*-2]}{2^*}$, then $W(x)\in L^1(\mathbb{R}^N)\cap \{L^q(\mathbb{R}^N)+L^{\infty}(\mathbb{R}^N)\}$ for some $q>\frac{2^*}{[\max(2\alpha, 1) 2^*-2]}$; Subcase (ii) $\alpha\geq \frac{N-1}{N}$ and $p<N$, then $W(x)\in L^1(\mathbb{R}^N)\cap \{L^q(\mathbb{R}^N)+L^{\infty}(\mathbb{R}^N)\}$ for some $q>1$. In any subcase, the solution is global existence for any initial data $u_0\in X$ by the results of Theorem 1. Moreover, if $0<\alpha<\frac{N-1}{N}$, then $p_c=\frac{N[\max(2\alpha,1)2^*-2]}{2^*}$ and
$q_c=\frac{2^*}{[\max(2\alpha,1)2^*-2]}$.\hfill$\Box$

{\bf Remark 4.1.} We would like to say something about the conditions on the initial $u_0$ if $W\in \{L^1(\mathbb{R}^N)\cap L^q(\mathbb{R}^N)\}$  in the critical case $q=q_c$. By the results of Theorem 1, if $a[h(s)+s^{\frac{1}{2}}]\geq \max(s^{\alpha},s^{\frac{1}{2}})$, $0<\alpha<\frac{N-1}{N}$, and $$a^22^{\frac{(q-1)N+2q}{qN}}C_s^{\frac{2}{2^*}}\|W\|_{L^q}\|u_0\|_2^{\min(4-4\alpha,2)}<1,$$ then
the solution is global existence when $q_c=\frac{2^*}{\max(2\alpha, 1) 2^*-2}$. On the other hand, in the case of $q_c=\frac{2^*}{\max(2\alpha, 1) 2^*-2}$, $E(u_0)\leq 0$ implies that
\begin{align*}
&\quad\int_{\mathbb{R}^N}|\nabla u_0|^2dx+\int_{\mathbb{R}^N}|\nabla h(|u_0|^2)|^2dx\leq \frac{1}{2}\int_{\mathbb{R}^N}(W*|u|^2)|u|^2dx\\
&\leq a^22^{\frac{(q-1)N+2q}{qN}}C_s^{\frac{2}{2^*}}\left(\int_{\mathbb{R}^N}|W|^qdx\right)^{\frac{1}{q}}
\left(\int_{\mathbb{R}^N}|u_0|^2dx\right)^{\frac{2q-1}{q\tau_1}}\\
&\qquad\times\left(\int_{\mathbb{R}^N}[|\nabla h(|u_0|^2)|^2+|\nabla u_0|^2]dx\right)^{\frac{2^*(2q-1)}{2q\tau_2}},
\end{align*}
which means that
$$
a^22^{\frac{(q-1)N+2q}{qN}}C_s^{\frac{2}{2^*}}\|W\|_{L^q}\|u_0\|_2^{\min(4-4\alpha,2)}\geq 1.
$$
Obviously, $S_{gl}(u_0)\cap S_{bl}(u_0)=\varnothing$, where
$$
S_{gl}(u_0)=\{u_0\in X,\ a^22^{\frac{((q-1)N+2q}{qN}}C_s^{\frac{2}{2^*}}\|W\|_{L^q}\|u_0\|_2^{\min(4-4\alpha,2)}<1\}
$$
and
\begin{align*}
S_{bl}(u_0)&=\{u_0\in X,\ E(u_0)\leq 0\}\\
&\subseteq\{u_0\in X,\ a^22^{\frac{(q-1)N+2q}{qN}}C_s^{\frac{2}{2^*}}\|W\|_{L^q}\|u_0\|_2^{\min(4-4\alpha,2)}\geq 1\}.
\end{align*}
That is, if $a[h(s)+s^{\frac{1}{2}}]\geq \max(s^{\alpha},s^{\frac{1}{2}})$ and $W(x)\in \{L^1(\mathbb{R}^N\cap L^q(\mathbb{R}^N)\}$ in the critical case of $q=q_c$,
\begin{align}
a^22^{\frac{(q-1)N+2q}{qN}}C_s^{\frac{2}{2^*}}\|W\|_{L^q}\|u_0\|_2^{\min(4-4\alpha,2)}=1\label{cztj1}
 \end{align}
 can be regarded as the watershed for the initial data $u_0$ which determines
whether the solution is global existence or not.

The similar conclusion for the initial data $u_0$ is also true for the following system, which is a special case of that in \cite{SW1}.
\begin{equation}
\label{cztj2} \left\{
\begin{array}{lll}
iu_t=\Delta u+2uh'(|u|^2) \Delta h(|u|^2)+F(|u|^2)u \quad {\rm for} \ x\in \mathbb{R}^N, \ t>0\\
u(x,0)=u_0(x),\quad x\in \mathbb{R}^N.
\end{array}\right.
\end{equation}
Assume that there exist positive constants $c_1$, $c_2$, $0<\theta<1$, $q>1$ satisfying the critical condition $(2-2^*)\theta+2q=2^*$ such that
$[G(s)]^{\theta}\leq c_1 s$ and $[G(s)]^q\leq c_2[s^{\frac{1}{2}}+h(s)]^{2^*}$ for $s\geq 0$. Here $c_1$ and $c_2$ are taken the values in the sense of the supremum.
Then the solution is global existence if $$2^{\frac{2^*-1}{\tilde{\tau}_2}}c_1^{\frac{1}{\tilde{\tau}_1}}(c_2C_s)^{\frac{1}{\tilde{\tau}_2}}\|u_0\|_2^{\frac{2}{\tilde{\tau}_1}}<1,$$
while the solution will blow up in finite time if $E(u_0)\leq 0$, $\Im \int_{\mathbb{R}^N}\bar{u}_0(x\cdot \nabla u_0)dx>0$, $|x|u_0\in L^2(\mathbb{R}^N)$
with other assumptions on $h(s)$ and $G(s)$. Here
$$
\frac{1}{\tilde{\tau}_1}=\frac{q-1}{q-\theta},\qquad \frac{1}{\tilde{\tau}_2}=\frac{1-\theta}{q-\theta}.
$$
But $E(u_0)\leq 0$ implies that
\begin{align*}
&\quad\int_{\mathbb{R}^N}|\nabla u_0|^2dx+\int_{\mathbb{R}^N}|\nabla h(|u_0|^2)|^2dx\leq \int_{\mathbb{R}^N}G(|u_0|^2)dx\\
&\leq \left(\int_{\mathbb{R}^N}c_1|u_0|^2dx\right)^{\frac{1}{\tilde{\tau}_1}}
\left(\int_{\mathbb{R}^N}c_2[|u_0|+h(|u_0|^2)]^{2^*}dx\right)^{\frac{1}{\tilde{\tau}_2}}\\
&\leq 2^{\frac{2^*-1}{\tilde{\tau}_2}}c_1^{\frac{1}{\tilde{\tau}_1}}(c_2C_s)^{\frac{1}{\tilde{\tau}_2}}\|u_0\|_2^{\frac{2}{\tilde{\tau}_1}}\int_{\mathbb{R}^N}[|\nabla u_0|^2+|\nabla h(|u_0|^2)|^2]dx
\end{align*}
which means that
$$
2^{\frac{2^*-1}{\tilde{\tau}_2}}c_1^{\frac{1}{\tilde{\tau}_1}}(c_2C_s)^{\frac{1}{\tilde{\tau}_2}}\|u_0\|_2^{\frac{2}{\tilde{\tau}_1}}\geq 1.
$$
Obviously, $S_{gl}(u_0)\cap S_{bl}(u_0)=\varnothing$, where
$$
S_{gl}(u_0)=\{u_0\in X,\quad 2^{\frac{2^*-1}{\tilde{\tau}_2}}c_1^{\frac{1}{\tilde{\tau}_1}}(c_2C_s)^{\frac{1}{\tilde{\tau}_2}}\|u_0\|_2^{\frac{2}{\tilde{\tau}_1}}<1\}
$$
and
$$
S_{bl}(u_0)=\{u_0\in X,\ E(u_0)\leq 0\}\subseteq \{u_0\in X,\quad 2^{\frac{2^*-1}{\tilde{\tau}_2}}c_1^{\frac{1}{\tilde{\tau}_1}}(c_2C_s)^{\frac{1}{\tilde{\tau}_2}}\|u_0\|_2^{\frac{2}{\tilde{\tau}_1}}\geq 1\}.
$$
That is,  in the critical case of $(2-2^*)\theta+2q=2^*$,
 \begin{align}
2^{\frac{2^*-1}{\tilde{\tau}_2}}c_1^{\frac{1}{\tilde{\tau}_1}}(c_2C_s)^{\frac{1}{\tilde{\tau}_2}}\|u_0\|_2^{\frac{2}{\tilde{\tau}_1}}=1\label{cztj3}
 \end{align}
 can be regarded as the watershed for the initial data $u_0$ which determines
whether the solution is global existence or not.

{\bf Remark 4.2.} If $h(|u|^2)=b|u|^{2\alpha}$, we would like to compare the results of \cite{SW1} with these of this paper when $W(x)\in L^1(\mathbb{R}^N)\cap \{L^q(\mathbb{R}^N)+L^{\infty}(\mathbb{R}^N)\}$ for some $q>1$, and discuss the similar roles of the exponents.
Consider
\begin{equation}
\label{mhs1} \left\{
\begin{array}{lll}
iu_t=\Delta u+2b^2\alpha |u|^{2\alpha-2}u \Delta (|u|^{2\alpha})+(|u|^{2\tilde{p}-2})u \quad {\rm for} \ x\in \mathbb{R}^N, \ t>0\\
u(x,0)=u_0(x),\quad x\in \mathbb{R}^N.
\end{array}\right.
\end{equation}

1. (i) If $b=0$, $h(s)\equiv 0$, the watershed exponent of $p$ for (\ref{mhs1}) in \cite{SW1} is $\tilde{p}_c=1+\frac{2}{N}$ in the means that the solution is global existence for any $u_0\in X$ if $\tilde{p}<\tilde{p}_c$ and the solution will blow up in finite time under certain conditions if $\tilde{p}> \tilde{p}_c$,
while the watershed exponent of $q$(or $p$) for (\ref{mhs2}) in this paper is $q_c=\frac{N}{2}$(or $p_c=2$)  in the means that the solution is global existence for any $u_0\in X$ if $q>q_c$(or $p<2$) and the solution will blow up in finite time under certain conditions if $q<q_c$(or $p>2$).

(ii) If $b>0$, $0<\alpha<\frac{N-1}{N}$ and $h(s)\not\equiv 0$, the watershed exponent of $p$ in (\ref{mhs1}) is $\tilde{p}_c=\max\{\frac{2}{N}, 2\alpha-1+\frac{2}{N}\}$ in the means that the solution is global existence for any $u_0\in X$ if $\tilde{p}<\tilde{p}_c$ and the solution will blow up in finite time under certain conditions if $\tilde{p}>\tilde{p}_c$, while the watershed exponent of $q$ for (\ref{mhs2}) is $q_c=\frac{2^*}{\max(2\alpha, 1) 2^*-2}$(or $p_c=\frac{N[\max(2\alpha,1)\cdot 2^*-2]}{2^*}$) in the means that the solution is global existence for any $u_0\in X$ if $q>q_c$(or $p<p_c$) and the solution will blow up in finite time under certain conditions if $q< q_c$(or $p>p_c$).

The watershed role of $q_c$(or $p_c$) for (\ref{mhs2}) in this paper is similar to that of $\tilde{p}_c$ for (\ref{mhs1}) in \cite{SW1}.

2. (i) If $b=0$, $h(s)\equiv 0$, the Sobolev critical exponent of $\tilde{p}$ for (\ref{mhs1}) is $\tilde{p}_s=\frac{2^*}{2}$, while $q_s=\frac{N}{4}$(or $p_s=4$) plays the similar role of Sobolev critical exponent for (\ref{mhs2}) when $N\geq 4$.
We can establish the blowup results on (\ref{mhs1}) if $\tilde{p}_c< \tilde{p}<\tilde{p}_s$ and those on (\ref{mhs2}) if $\frac{N}{4}=q_s<q<q_c=\frac{N}{2}$(or $2=p_c<p<p_s=4$) for some initial $u_0$.

(ii) If $b>0$, $0<\alpha<\frac{N-1}{N}$ and $h(s)\not\equiv 0$, the Sobolev critical exponent of $\tilde{p}$ for (\ref{mhs1}) is $\tilde{p}_s=\max(\alpha, \frac{1}{2})2^*$, while $q_s=\max(\frac{\max(\alpha,\frac{1}{2})2^*}{[\max(2\alpha,1)\cdot 2^*-2]},1)$(or $p_s=N\min(\frac{[\max(2\alpha,1)\cdot 2^*-2]}{\max(\alpha,\frac{1}{2})2^*},1)$) likes Sobolev critical exponent for (\ref{mhs2}).
We can establish the blowup results on (\ref{mhs1}) if $\tilde{p}_c< \tilde{p}<\tilde{p}_s$ and those on (\ref{mhs2}) if $\max(\frac{\max(\alpha,\frac{1}{2})2^*}{[\max(2\alpha,1)\cdot 2^*-2]},1)=q_s<q< q_c=\frac{2^*}{[\max(2\alpha,1)\cdot 2^*-2]}$(or $\frac{N[\max(2\alpha,1)\cdot 2^*-2]}{2^*}=p_c< p<p_s=N(\min(\frac{[\max(2\alpha,1)\cdot 2^*-2]}{\max(\alpha,\frac{1}{2})2^*},1)$) for some initial $u_0$.

(iii) We also point out that $\tilde{p}_c<\tilde{p}_s$, $p_c<p_s$, while $q_s<q_c$. The Sobolev critical exponent role of $q_s$(or $p_s$) for (\ref{mhs2}) in this paper is similar to that of $\tilde{p}_s$ for (\ref{mhs1}) in \cite{SW1}.

\section{Sharp threshold for the global existence and blowup in finite time for the solution of (\ref{1})}

\qquad In this section, we will  establish a sharp threshold for the blowup and global existence of the solution to (\ref{1}) under certain conditions.

{\bf Theorem 3.} (Sharp Threshold ) {\it Let $u(x,t)$ be the solution of (\ref{1}) with $u_0\in X$.
 Assume that:

 (i) There exists constant $k\in \mathbb{R}$ such that
$sh''(s)\leq kh'(s)$ if $h'(s)\geq 0$ or $sh''(s)\geq kh'(s)$ if $h'(s)\leq 0$, there exist $a>0$ in the sense of infimum and $0<\alpha<\frac{N-1}{N}$ in the sense of supremum such that $\max(s^{\alpha},s^{\frac{1}{2}})\leq a[h(s)+s^{\frac{1}{2}}]$ for $s\geq 0$, and
$$(2-\underline{l})+4(N+2-\underline{l})[h'(s)]^2s+8Nh''(s)h'(s)s^2\geq 0$$
for $0<\underline{l}\leq 2$.

(ii) $W(x)\geq 0$ for $x\in \mathbb{R}^N$, $W(x)\in \{L^1(\mathbb{R}^N\cap L^q(\mathbb{R}^N)\}$, $\max(1,\frac{N}{4})<q<\frac{N}{2}$ if $h(s)\equiv0$, $\max\left(1,\frac{\max(\alpha,\frac{1}{2})2^*}{[\max(2\alpha,1)2^*-2]}\right)<q<\frac{2^*}{[\max(2\alpha,1)2^*-2]}$ if $h(s)\not\equiv 0$ and $0<\alpha<\frac{N-1}{N}$, and there exist constant $L>1+\frac{\max[(2k+1)N,0]}{2}$ and $C$
such that
$$ LW(x)\leq -\frac{x\cdot \nabla W}{2}\leq CW(x).$$

Moreover, suppose that
there exists $\omega>0$ such that
\begin{align} d_I:=\inf_{\{w\in H^1(\mathbb{R}^N)\setminus \{0\};
Q(w)=0\}}\left(\frac{\omega}{2}\|w\|_2^2+E(w)\right)>0,\label{9651}\end{align}
where \begin{align}
Q(w)&=2\int_{\mathbb{R}^N}|\nabla w|^2dx+(N+2)\int_{\mathbb{R}^N}|\nabla h(|w|^2)|^2dx\nonumber\\
&\quad +8N\int_{\mathbb{R}^N}h''(|w|^2)h'(|w|^2)|w|^4|\nabla w|^2dx-\frac{1}{2}\int_{\mathbb{R}^N}[(x\cdot \nabla W)*|w|^2]|w|^2dx,\\
E(w)&=\frac{1}{2}\int_{\mathbb{R}^N}[|\nabla w|^2+|\nabla h(|w|^2)|^2]dx-\frac{1}{4}\int_{\mathbb{R}^N}(W*|w|^2)|w|^2dx,
\end{align}
and $u_0$ satisfies
$$\frac{\omega}{2}\|u_0\|_2^2+E(u_0)<d_I.$$

Then we have:

(1). If $Q(u_0)>0$, the solution of (\ref{1}) exists
globally;

(2). If $Q(u_0)<0$ and $\Im \int_{\mathbb{R}^N} \bar{u}_0(x\cdot \nabla u_0)dx\geq 0$, $xu_0\in L^2(\mathbb{R}^N)$, the solution
of (\ref{1}) blows up in finite time.}

{\bf The proof of Theorem 3. } We proceed in four Steps.

Step 1. We will prove $d_I>0$.

Since $Q(w)=0$, $w\not\equiv 0$, we have
\begin{align}
&\quad \underline{l}(\int_{\mathbb{R}^N}|\nabla w|^2dx+\int_{\mathbb{R}^N}|\nabla h(|w|^2)|^2dx)\nonumber\\
&=2\int_{\mathbb{R}^N}|\nabla w|^2dx+(N+2)\int_{\mathbb{R}^N}|\nabla h(|w|^2)|^2dx+8N\int_{\mathbb{R}^N}h''(|w|^2)h'(|w|^2)|w|^4|\nabla w|^2dx\nonumber\\
&\quad-\int_{\mathbb{R}^N}\left[(2-\underline{l})+4(N+2-\underline{l})(h'(|w|^2))^2|w|^2+8Nh''(|w|^2)h'(|w|^2)|w|^4\right]|\nabla w|^2dx\nonumber\\
&\leq -\frac{1}{2}\int_{\mathbb{R}^N}[(x\cdot \nabla W)*|w|^2]|u|^2dx\leq C\int_{\mathbb{R}^N}(W*|w|^2)|w|^2dx\nonumber\\
&\leq C\left(\int_{\mathbb{R}^N}|W|^qdx\right)^{\frac{1}{q}}\left(\int_{\mathbb{R}^N}|w|^{\frac{4q}{2q-1}}\right)^{\frac{2q-1}{q}}\nonumber\\
&\leq C\left[\left(\int_{\mathbb{R}^N}|w|^2dx\right)^{\frac{1}{\tau_1}}\left(\int_{\mathbb{R}^N}|w|^{\max(2\alpha, 1)2^*}dx\right)^{\frac{1}{\tau_2}}\right]^{\frac{2q-1}{q}}\nonumber\\
&\leq C\left(\int_{\mathbb{R}^N}|w|^2dx\right)^{\frac{2q-1}{q\tau_1}}\left(\int_{\mathbb{R}^N}a^{2^*}[|w|^{2^*}+|h(|w|^2)|^{2^*}]dx\right)^{\frac{2q-1}{q\tau_2}}\nonumber\\
&\leq C'\left(\int_{\mathbb{R}^N}|w|^2dx\right)^{\frac{2q-1}{q\tau_1}}\left(\int_{\mathbb{R}^N}|\nabla w|^2dx+\int_{\mathbb{R}^N}|\nabla h(|w|^2)|^2dx\right)^{\frac{2^*(2q-1)}{2q\tau_2}},\label{6291}
\end{align}
which implies that
$$
\left(\int_{\mathbb{R}^N}|w|^2dx\right)^{\frac{2q-1}{q\tau_1}}\left(\int_{\mathbb{R}^N}|\nabla w|^2dx+\int_{\mathbb{R}^N}|\nabla h(|w|^2)|^2dx\right)^{\frac{2^*(2q-1)}{2q\tau_2}-1}\geq C>0.
$$
Since $q< \frac{2^*}{[\max(2\alpha,1)2^*-2]}$ equals to $\frac{2^*(2q-1)}{2q\tau_2}-1>0$, we have
\begin{align}
\int_{\mathbb{R}^N}|w|^2dx+\int_{\mathbb{R}^N}|\nabla w|^2dx+\int_{\mathbb{R}^N}|\nabla h(|w|^2)|^2dx\geq C>0.\label{1018w3}
\end{align}
Here
$$
 \frac{1}{\tau_1}=\frac{(2q-1)[\max(2\alpha,1)2^*-2]-2}{(2q-1)[\max(2\alpha,1)2^*-2]},\qquad \frac{1}{\tau_2}=\frac{2}{(2q-1)[\max(2\alpha,1)2^*-2]}.
$$

On the other hand, using $Q(w)=0$ again, we get
\begin{align}
&\quad(\max[(2k+1)N,0]+2)\int_{\mathbb{R}^N}[|\nabla w|^2+|\nabla h(|w|^2)|^2]dx\nonumber\\
&\geq 2\int_{\mathbb{R}^N}[|\nabla w|^2+((2k+1)N+2)\int_{\mathbb{R}^N}|\nabla h(|w|^2)|^2]dx\nonumber\\
&\geq 2\int_{\mathbb{R}^N}|\nabla w|^2dx+(N+2)\int_{\mathbb{R}^N}|\nabla h(|w|^2)|^2dx+8N\int_{\mathbb{R}^N}h''(|w|^2)h'(|w|^2)|w|^4|\nabla w|^2dx\nonumber\\
&=-\frac{1}{2}\int_{\mathbb{R}^N}[(x\cdot \nabla W)*|u|^2]|u|^2dx\geq L\int_{\mathbb{R}^N}(W*|u|^2)|u|^2dx.\label{1019s1}
\end{align}
Therefore
\begin{align}
E(w)&=\frac{1}{2}\int_{\mathbb{R}^N}[|\nabla w|^2dx+\int_{\mathbb{R}^N}|\nabla h(|w|^2)|^2]dx-\frac{1}{4}\int_{\mathbb{R}^N}(W*|w|^2)|w|^2dx\nonumber\\
&\geq \frac{1}{2}\left(1-\frac{(\max[(2k+1)N,0]+2)}{2L}\right)\int_{\mathbb{R}^N}[|\nabla w|^2+|\nabla h(|w|^2)|^2]dx.\label{1018w4}
\end{align}
(\ref{1019s1}) and (\ref{1018w4}) mean that
\begin{align}
&\quad\frac{\omega}{2} \int_{\mathbb{R}^N}|w|^2dx+E(w)\nonumber\\
&\geq \frac{1}{2}\min\left(\omega, 1-\frac{(\max[(2k+1)N,0]+2)}{2L}\right)\int_{\mathbb{R}^N}[|w|^2+|\nabla w|^2+|\nabla h(|w|^2)|^2]dx\nonumber\\
&\geq C>0.\label{1018w3'}
\end{align}
Therefore $d_I>0$.

Step 2. Denote
$$
K_+=\{u\in H^1(\mathbb{R}^N)\setminus\{0\},\ Q(u)>0,\
\frac{\omega}{2}\|u\|_2^2+E(u)<d_I\}
$$
and
$$
K_-=\{u\in H^1(\mathbb{R}^N)\setminus\{0\},\ Q(u)<0,\
\frac{\omega}{2}\|u\|_2^2+E(u)<d_I\}.
$$
We will prove that $K_+$ and $K_-$ are invariant sets of
(\ref{1}).

Assume that $u_0\in K_+$, i.e., $Q(u_0)>0$ and
$\frac{\omega}{2}\|u_0\|_2^2+E(u_0)<d_I$. It is easy to verify that
\begin{align}
u(\cdot,t)\in H^1(\mathbb{R}^N)\setminus\{0\}, \quad
\frac{\omega}{2}\|u(\cdot,t)\|_2^2+E(u(\cdot,t))<d_I.\label{9626z2}\end{align}
because $\|u\|^2_2$ and $E(u)$ are conservation quantities for (\ref{1}).

We need to show that $Q(u(\cdot,t))>0$ for $t\in (0,T)$. Contradictorily, if there exists $t_1\in (0,T)$ such that
$Q(u(\cdot,t_1))<0$, then there exists a $t_2\in [0, t_1]$ such that
$Q(u(\cdot,t_2))=0$ by the continuity. And
$$\frac{\omega}{2}\|u(\cdot,t_2)\|_2^2+E(u(\cdot,t_2))<d_I$$
by (\ref{9626z2}), which is a contradiction to the definition of $d_I$.
Hence $Q(u(\cdot,t))>0$. This inequality and (\ref{9626z2}) imply
that $u(\cdot,t)\in K_+$, which means that $K_+$ is a invariant
set of (\ref{1}).

Similarly, we can prove that $K_-$ is also a
invariant set of (\ref{1}). We omit the details here.

Step 3. Assume that $Q(u_0)>0$ and
$\frac{\omega}{2}\|u_0\|_2^2+E(u_0)<d_I$. Since $\mathcal{K}$ is invariant set of (\ref{1}), we
have $Q(u(\cdot,t))>0$ and
$\frac{\omega}{2}\|u(\cdot,t)\|_2^2+E(u(\cdot,t))<d_I$.  Using $Q(u(\cdot,t))>0$, we get
\begin{align}
&\quad(\max[(2k+1)N,0]+2)\int_{\mathbb{R}^N}[|\nabla u(\cdot,t)|^2+|\nabla h(|u(\cdot,t)|^2)|^2]dx\nonumber\\
&\geq 2\int_{\mathbb{R}^N}[|\nabla u(\cdot,t)|^2+((2k+1)N+2)\int_{\mathbb{R}^N}|\nabla h(|u(\cdot,t)|^2)|^2]dx\nonumber\\
&\geq 2\int_{\mathbb{R}^N}|\nabla u(\cdot,t)|^2dx+(N+2)\int_{\mathbb{R}^N}|\nabla h(|u(\cdot,t)|^2)|^2dx\nonumber\\
&\quad+8N\int_{\mathbb{R}^N}h''(|u(\cdot,t)|^2)h'(|u(\cdot,t)|^2)|u(\cdot,t)|^4|\nabla u(\cdot,t)|^2dx\nonumber\\
&=-\frac{1}{2}\int_{\mathbb{R}^N}[(x\cdot \nabla W)*|u(\cdot,t)|^2]|u(\cdot,t)|^2dx\geq L\int_{\mathbb{R}^N}(W*|u(\cdot,t)|^2)|u(\cdot,t)|^2dx.\label{1019s1}
\end{align}
Therefore
\begin{align}
E(u(\cdot,t))&=\frac{1}{2}\int_{\mathbb{R}^N}[|\nabla u(\cdot,t)|^2dx+|\nabla h(|u(\cdot,t)|^2)|^2]dx-\frac{1}{4}\int_{\mathbb{R}^N}(W*|u(\cdot,t)|^2)|u(\cdot,t)|^2dx\nonumber\\
&\geq \frac{2L-(\max[(2k+1)N,0]+2)}{4L}\int_{\mathbb{R}^N}[|\nabla u(\cdot,t)|^2+|\nabla h(|u(\cdot,t)|^2)|^2]dx.\label{1018w4}
\end{align}
By the mass and energy conservation laws, (\ref{1019s1}) and (\ref{1018w4}) mean that
\begin{align*}
&d_I\geq \frac{\omega}{2} \int_{\mathbb{R}^N}|u_0|^2dx+E(u_0)=\frac{\omega}{2} \int_{\mathbb{R}^N}|u(\cdot,t)|^2dx+E(u(\cdot,t))\nonumber\\
&\geq \frac{1}{2}\min\left(\omega, 1-\frac{(\max[(2k+1)N,0]+2)}{2L}\right)\int_{\mathbb{R}^N}[|u(\cdot,t)|^2+|\nabla u(\cdot,t)|^2+|\nabla h(|u(\cdot,t)|^2)|^2]dx\nonumber\\
&\geq C>0.
\end{align*}
and
$$
\int_{\mathbb{R}^N}|\nabla u(\cdot,t)|^2dx+\int_{\mathbb{R}^N}|\nabla h(|u(\cdot,t)|^2)|^2dx\leq C<\infty,
$$
i.e., the solution $u(x,t)$ of (\ref{1}) exists globally.

Step 4. Suppose that $|x|u_0\in L^2(\mathbb{R}^N)$, $Q(u_0)<0$ and
$\frac{\omega}{2}\|u_0\|_2^2+E(u_0)<d_I$. Since $\mathcal{K}_-$ is a invariant set of (\ref{1}), we have $Q(u(\cdot,t))<0$ and
$\frac{\omega}{2}\|u(\cdot,t)\|_2^2+E(u(\cdot,t))<d_I$.

Let $J(t)=\int_{\mathbb{R}^N} |x|^2|u|^2dx$. Then
$$
J''(t)=4Q(u(x,t)),\quad J'(t)=-4\Im \int_{\mathbb{R}^N} \bar{u}(x\cdot \nabla u)dx.
$$
Since
$J'(0)=-4\Im \int_{\mathbb{R}^N} \bar{u}_0(x\cdot \nabla u_0)dx<0$, we have
$$
J'(t)=J'(0)+\int_0^t J''(\tau)d\tau=J'(0)+4\int_0^t Q(u(\cdot,\tau))d\tau<J'(0)<0
$$
and
$$
0\leq J(t)=J(0)+\int_0^t J'(\tau)d\tau<J(0)+J'(0)t,
$$
which implies that the maximum existence interval for $t$ is finite, and the solution blows up in finite time.\hfill
$\Box$

{\bf Remark 5.1.} We give an example below. If $h(s)=s^{\alpha}$, $\frac{1}{2}<\alpha$ and
closes to $\frac{1}{2}$ enough, then $k=\alpha-1$. Let
\begin{equation*}
W(x)=W(r)=\left\{
\begin{array}{lll}
\frac{1}{r^{[(2\alpha-1)N+2+\epsilon]}},\quad 0<r\leq 1\\
f(r),\quad 1\leq r\leq 2\\
\frac{1}{r^{C[(2\alpha-1)N+2+\epsilon]}},\quad r\geq 2.
\end{array}\right.
\end{equation*}
Here $\epsilon$ is a positive constant small enough, $f(r)$ is an arbitrary function satisfying $f(1)=1$, $f(2)=\frac{1}{2^{C[(2\alpha-1)N+2+\epsilon]}}$ and
$$
2Lf(r)+rf'(r)\leq 0\leq 2Cf(r)+rf'(r), \quad 1\leq r\leq 2.
$$
Then $W(x)\in \{L^1(\mathbb{R}^N\cap L^q(\mathbb{R}^N)\}$ for some $\max\left(1,\frac{\max(\alpha,\frac{1}{2})2^*}{[\max(2\alpha,1)2^*-2]}\right)<q<\frac{2^*}{[\max(2\alpha,1)2^*-2]}$.
The assumptions on $h(s)$ and $W(x)$ can be satisfied for suitable $L$ and $C$, and the conclusions of Theorem 3 are true under certain conditions on initial data.

\section{The pseudo-conformal conservation laws and asymptotic behavior for the solution}
\qquad In this section, we give some extensions of the work in \cite{Ginibre1, Ginibre2} to the quasilinear case, establishing two pseudo-conformal conservation laws, which in turn yield asymptotic behavior for the solution of (\ref{1}).

{\bf Theorem 4.( Pseudo-conformal Conservation Laws)} {\it 1. Assume that $u$ is the global solution of (\ref{1}), $u_0\in X$ and $xu_0\in L^2(\mathbb{R}^N)$.  Then
\begin{align}
P(t)&=\int_{\mathbb{R}^N}|(x-2it\nabla)u|^2dx+4t^2\int_{\mathbb{R}^N}|\nabla h(|u|^2)|^2dx-2t^2\int_{\mathbb{R}^N}(W*|u|^2)|u|^2dx\nonumber\\
&=\int_{\mathbb{R}^N}|xu_0|^2dx+4\int_0^t\tau\theta(\tau)d\tau.\label{691}
\end{align}

2.  Assume that $u$ is the blowup solution of (\ref{1}) with blowup time $T$, $u_0\in X$ and $xu_0\in L^2(\mathbb{R}^N)$. Then

\begin{align}
B(t)&:=\int_{\mathbb{R}^N}|(x+2i(T-t)\nabla)u|^2dx+4(T-t)^2\int_{\mathbb{R}^N}|\nabla h(|u|^2)|^2dx\nonumber\\
&\quad-2(T-t)^2\int_{\mathbb{R}^N}(W*|u|^2)|u|^2dx\nonumber\\
&=\int_{\mathbb{R}^N}|(x+2iT\nabla)u_0|^2dx+4T^2\int_{\mathbb{R}^N}|\nabla h(|u_0|^2)|^2dx-2T^2\int_{\mathbb{R}^N}(W*|u_0|^2)|u_0|^2dx\nonumber\\
&\quad -4\int_0^t(T-\tau)\theta(\tau)d\tau.\label{893}
\end{align}

Here
\begin{align}
\theta(t)&=\int_{\mathbb{R}^N}-4N[2h''(|u|^2)h'(|u|^2)|u|^2+( h'(|u|^2))^2]|u|^2|\nabla u|^2dx\nonumber\\
&\qquad-\int_{\mathbb{R}^N} \left([W+\frac{(x\cdot \nabla W)}{2}]*|u|^2\right)|u|^2dx.\label{691'}
\end{align}

}

{\bf Proof of Theorem 4:} 1. Assume that $u$ is the global solution of (\ref{1}), $u_0\in X$ and $xu_0\in L^2(\mathbb{R}^N)$. Using the conservation of energy,
we have
\begin{align}
P(t)&=\int_{\mathbb{R}^N}|xu|^2dx+4t\Im \int_{\mathbb{R}^N}\bar{u}(x\cdot \nabla u)dx+4t^2\int_{\mathbb{R}^N}|\nabla u|^2dx\nonumber\\
&\qquad+4t^2\int_{\mathbb{R}^N}|\nabla h(|u|^2)|^2dx-2t^2\int_{\mathbb{R}^N}(W*|u|^2)|u|^2dx\nonumber\\
&=\int_{\mathbb{R}^N}|xu|^2dx+4t\Im \int_{\mathbb{R}^N}\bar{u}(x\cdot \nabla u)dx+8t^2E(u_0).\label{692}
\end{align}
Noticing that
$$\frac{d}{dt} \int_{\mathbb{R}^N}|x|^2|u|^2dx=-4\Im \int_{\mathbb{R}^N} \bar{u}(x\cdot \nabla u)dx,$$
we obtain
\begin{align}
P'(t)&=\frac{d}{dt}\int_{\mathbb{R}^N}|xu|^2dx+4\Im \int_{\mathbb{R}^N}\bar{u}(x\cdot \nabla u)dx+4t\frac{d}{dt}\Im \int_{\mathbb{R}^N}\bar{u}(x\cdot \nabla u)dx+16tE(u_0)\nonumber\\
&=4t\frac{d}{dt}\Im \int_{\mathbb{R}^N}\bar{u}(x\cdot \nabla u)dx+16tE(u_0)\nonumber\\
&=4t\left\{-2\int_{\mathbb{R}^N}|\nabla u|^2dx-(N+2)\int_{\mathbb{R}^N}|\nabla h(|u|^2)|^2dx\right.\nonumber\\
&\qquad \left.-8N\int_{\mathbb{R}^N}h''(|u|^2)h'(|u|^2)|u|^4|\nabla u|^2dx-\frac{1}{2}\int_{\mathbb{R}^N}[(x\cdot \nabla W)*|u|^2]|u|^2dx\right\}\nonumber\\
&\qquad+8t\int_{\mathbb{R}^N}[|\nabla u|^2+|\nabla h(|u|^2)|^2]dx-4t\int_{\mathbb{R}^N}(W*|u|^2)|u|^2dx\nonumber\\
&=4t\int_{\mathbb{R}^N}-4N[2h''(|u|^2)h'(|u|^2)|u|^2+(h'(|u|^2))^2]|u|^2|\nabla u|^2dx\nonumber\\
&\quad-4t\int_{\mathbb{R}^N}([W+\frac{(x\cdot \nabla W)}{2}]*|u|^2)|u|^2dx\nonumber\\
&=4t\theta(t).\label{693}
\end{align}
Integrating (\ref{693}) from $0$ to $t$, we have
\begin{align}
P(t)&=\int_{\mathbb{R}^N}|(x-2it\nabla)u|^2dx+4t^2\int_{\mathbb{R}^N}|\nabla h(|u|^2)|^2dx-2t^2\int_{\mathbb{R}^N}(W*|u|^2)|u|^2dx\nonumber\\
&=\int_{\mathbb{R}^N}|xu_0|^2dx+4\int_0^t\tau\theta(\tau)d\tau,\label{694}
\end{align}
where $\theta(\tau)$ is defined by (\ref{691'}).

2. Assume that $u$ is the blowup solution of (\ref{1}),  $u_0\in X$ and $xu_0\in L^2(\mathbb{R}^N)$. By the conservation of energy, we have
\begin{align}
B(t)&:=\int_{\mathbb{R}^N}|(x+2i(T-t)\nabla)u|^2dx+4(T-t)^2\int_{\mathbb{R}^N}|\nabla h(|u|^2)|^2dx\nonumber\\
&\quad-2(T-t)^2\int_{\mathbb{R}^N}(W*|u|^2)|u|^2dx\nonumber\\
&=\int_{\mathbb{R}^N}|xu|^2dx-4(T-t)\Im \int_{\mathbb{R}^N}\bar{u}(x\cdot \nabla u)dx+4(T-t)^2\int_{\mathbb{R}^N}|\nabla u|^2dx\nonumber\\
&\qquad+4(T-t)^2\int_{\mathbb{R}^N}|\nabla h(|u|^2)|^2dx-2(T-t)^2\int_{\mathbb{R}^N}(W*|u|^2)|u|^2dx\nonumber\\
&=\int_{\mathbb{R}^N}|xu|^2dx-4(T-t)\Im \int_{\mathbb{R}^N}\bar{u}(x\cdot \nabla u)dx+8(T-t)^2E(u_0)\label{891}
\end{align}
and
\begin{align}
B'(t)&=\frac{d}{dt}\int_{\mathbb{R}^N}|xu|^2dx+4\Im \int_{\mathbb{R}^N}\bar{u}(x\cdot \nabla u)dx\nonumber\\
&\quad-4(T-t)\frac{d}{dt}\Im \int_{\mathbb{R}^N}\bar{u}(x\cdot \nabla u)dx-16(T-t)E(u_0)\nonumber\\
&=-4(T-t)\frac{d}{dt}\Im \int_{\mathbb{R}^N}\bar{u}(x\cdot \nabla u)dx-16(T-t)E(u_0)\nonumber\\
&=4(T-t)\left\{\int_{\mathbb{R}^N}4N[2h''(|u|^2)h'(|u|^2)|u|^2+( h'(|u|^2))^2]|u|^2|\nabla u|^2dx\right.\nonumber\\
&\qquad\left.+\int_{\mathbb{R}^N} ([W+\frac{(x\cdot \nabla W)}{2}]*|u|^2)|u|^2dx\right\}.\label{892}
\end{align}
Integrating (\ref{892}) from $0$ to $t$, we have
\begin{align*}
B(t)&=B(0)-4\int_0^t(T-\tau)\theta(\tau)d\tau\nonumber\\
&=\int_{\mathbb{R}^N}|(x+2iT\nabla)u_0|^2dx+4T^2\int_{\mathbb{R}^N}[|\nabla h(|u_0|^2)|^2dx\nonumber\\
&\quad-2T^2\int_{\mathbb{R}^N}(W*|u_0|^2)|u_0|^2dx-4\int_0^t(T-\tau)\theta(\tau)d\tau,
\end{align*}
where $\theta(\tau)$ is defined by (\ref{691'}).\hfill $\Box$

As the applications of Theorem 4, we give some asymptotic behavior results on the global solution of (\ref{1}) and the lower bound for the blowup rate the blowup solution of  (\ref{1}).

{\bf Theorem 5.} {\it 1. Assume that $u$ is the global solution of (\ref{1}),  $u_0\in X$, $xu_0\in L^2(\mathbb{R}^N)$, and $W(x)\leq 0$ for $x\in \mathbb{R}^N$. Then the following properties hold:

(1) If $2h''(s)h'(s)s+(h'(s))^2\geq 0$ for $s\geq 0$, and $2W+(x\cdot \nabla W)\geq 0$ for $x\in \mathbb{R}^N$, then there exists $C$ such that
\begin{align}
\int_{\mathbb{R}^N}|\nabla h(|u|^2)|^2dx +\int_{\mathbb{R}^N}(|W|*|u|^2)|u|^2dx\leq Ct^{-2}\quad {\rm for}\ t\geq 1. \label{6101}
\end{align}

(2) If $2h''(s)h'(s)s+(h'(s))^2\geq 0$ for $s\geq 0$, and $-c|W|\leq 2W+(x\cdot \nabla W)\leq 0$ for $x\in \mathbb{R}^N$ for some $0<c<2$, then there exists $C$ such that
\begin{align}
\int_{\mathbb{R}^N}|\nabla h(|u|^2)|^2dx +\int_{\mathbb{R}^N}(|W|*|u|^2)|u|^2dx\leq\frac{C}{t^{2-c}}\quad {\rm for}\ t\geq 1.\label{6102}
\end{align}

(3) If $-k_1(h'(s))^2\leq 2h''(s)h'(s)s+(h'(s))^2\leq 0$ for some $0<k_1<\frac{2}{N}$, and $2W+(x\cdot \nabla W)\geq 0$ for $x\in \mathbb{R}^N$, then there exists $C$ such that
\begin{align}
\int_{\mathbb{R}^N}|\nabla h(|u|^2)|^2dx +\int_{\mathbb{R}^N}(|W|*|u|^2)|u|^2dx\leq \frac{C}{t^{2-Nk_1}}\quad {\rm for}\ t\geq 1.\label{6103}
\end{align}

(4) If $-k_1(h'(s))^2\leq 2h''(s)h'(s)s+(h'(s))^2\leq 0$ for some $0<k_1<\frac{2}{N}$, and $-c|W|\leq 2W+(x\cdot \nabla W)\leq 0$ for $x\in \mathbb{R}^N$ for some $0<c<2$, then there exists $C$ such that
\begin{align}
\int_{\mathbb{R}^N}|\nabla h(|u|^2)|^2dx +\int_{\mathbb{R}^N}(|W|*|u|^2)|u|^2dx\leq \frac{C}{t^{2-\max(Nk_1,c)}}\quad {\rm for}\ t\geq 1.\label{6103}
\end{align}

In all cases above, by the conservation of energy, we have
\begin{align}
\lim_{t\rightarrow +\infty} \int_{\mathbb{R}^N}|\nabla u(x,t)|^2dx=2E(u_0),\quad \lim_{t\rightarrow +\infty} \|u(\cdot,t)\|_{H^1}^2=m^2(u_0)+2E(u_0).\label{9181}
\end{align}

2. Assume that $u$ is the blowup solution of (\ref{1}), $[(h'(s))^2+2h''(s)h'(s)s]\leq 0$ for $s\geq 0$, $W(x)\geq 0$ and $2W+(x\cdot \nabla W)\leq 0$ for $x\in \mathbb{R}^N$, $u_0\in X$ and $xu_0\in L^2(\mathbb{R}^N)$. If $E(u_0)\leq 0$ and $-4T^2E(u_0)-\int_{\mathbb{R}^N}|xu_0|^2dx-4T\Im\int_{\mathbb{R}^N}\bar{u}_0(x\cdot \nabla u_0)dx>0$, then
\begin{align}
\int_{\mathbb{R}^N}[|\nabla u|^2+|\nabla h(|u|^2)|^2]dx\geq \frac{C}{(T-t)^2},\quad \int_{\mathbb{R}^N}(W*|u|^2)|u|^2dx\geq \frac{C}{(T-t)^2}.\label{895}
\end{align}

}

{\bf Proof of Theorem 5:} 1. Assume that $u$ is the global solution of (\ref{1}), $u_0\in X$ and $xu_0\in L^2(\mathbb{R}^N)$, $W(x)\leq 0$.

(1) $2h''(s)h'(s)s+(h'(s))^2\geq 0$ and $2W+(x\cdot \nabla W)\geq 0$. (\ref{691}) implies that
$$
4t^2\int_{\mathbb{R}^N}|\nabla h(|u|^2)|^2dx+2t^2\int_{\mathbb{R}^N}(|W|*|u|^2)|u|^2dx\leq\int_{\mathbb{R}^N}|xu_0|^2dx,
$$
i.e.,
$$
\int_{\mathbb{R}^N}|\nabla h(|u|^2)|^2dx+\int_{\mathbb{R}^N}(|W|*|u|^2)|u|^2dx\leq Ct^{-2}.
$$

(2) $2h''(s)h'(s)s+(h'(s))^2\geq 0$ and $-c|W|\leq 2W+(x\cdot \nabla W)\leq 0$ for some $0<c<2$. (\ref{691}) implies that
\begin{align}
&\quad 4t^2\int_{\mathbb{R}^N}|\nabla h(|u|^2)|^2dx+2t^2\int_{\mathbb{R}^N}(|W|*|u|^2)|u|^2dx\nonumber\\
&\leq\int_{\mathbb{R}^N}|xu_0|^2dx+2c\int_0^t\tau\left(\int_{\mathbb{R}^N}(|W|*|u|^2)|u|^2dx\right)d\tau.
\label{6232}
\end{align}
Let
$$
A_1(t):=2\int_0^t\tau\left(\int_{\mathbb{R}^N}(|W|*|u|^2)|u|^2dx\right)d\tau.
$$
(\ref{6232}) implies
$$
A'_1(t)\leq \frac{C_0}{t}+\frac{c}{t}A_1(t).
$$
Using Gronwall's inequality, we have
$$
A_1(t)\leq t^c[A_1(1)+C-\frac{C}{t^c}]\leq C't^c.
$$
Consequently,
$$
\int_{\mathbb{R}^N}|\nabla h(|u|^2)|^2dx+\int_{\mathbb{R}^N}(|W|*|u|^2)|u|^2dx\leq \frac{C}{t^{2-c}}.
$$

(3) $-k_1(h'(s))^2<2h''(s)h'(s)s+(h'(s))^2<0$ for some $0<k_1<\frac{2}{N}$ and $2W+(x\cdot \nabla W)\geq 0$.
(\ref{691}) implies that
\begin{align}
&\quad 4t^2\int_{\mathbb{R}^N}|\nabla h(|u|^2)|^2dx+2t^2\int_{\mathbb{R}^N}(|W|*|u|^2)|u|^2dx\nonumber\\
&\leq\int_{\mathbb{R}^N}|xu_0|^2dx+4Nk_1\int_0^t\tau\left(\int_{\mathbb{R}^N}|\nabla h(|u|^2)|^2dx\right)d\tau.
\label{914w1}
\end{align}
Let
$$
A_2(t):=4\int_0^t\tau\left(\int_{\mathbb{R}^N}|\nabla h(|u|^2)|^2dx\right)d\tau.
$$
(\ref{914w1}) implies
$$
A'_2(t)\leq \frac{C_0}{t}+\frac{Nk_1}{t}A_2(t).
$$
Using Gronwall's inequality, we have
$$
A_2(t)\leq t^{Nk_1}c[A_2(1)+C-\frac{C}{t^{Nk_1}}]\leq C't^{Nk_1}.
$$
Consequently, we have
$$
\int_{\mathbb{R}^N}|\nabla h(|u|^2)|^2dx+\int_{\mathbb{R}^N}(|W|*|u|^2)|u|^2dx\leq \frac{C}{t^{2-Nk_1}}.
$$

(4) $-k_1(h'(s))^2<2h''(s)h'(s)s+(h'(s))^2<0$ for some $0<k_1<\frac{2}{N}$ and $-c|W|\leq 2W+(x\cdot \nabla W)\leq 0$ for some $0<c<2$.
(\ref{691}) implies that
\begin{align}
&\quad 4t^2\int_{\mathbb{R}^N}|\nabla h(|u|^2)|^2dx+2t^2\int_{\mathbb{R}^N}(|W|*|u|^2)|u|^2dx\nonumber\\
&\leq\int_{\mathbb{R}^N}|xu_0|^2dx+4Nk_1\int_0^t\tau\left(\int_{\mathbb{R}^N}|\nabla h(|u|^2)|^2dx\right)d\tau
+2c\int_0^t\tau\left(\int_{\mathbb{R}^N}(|W|*|u|^2)|u|^2dx\right)d\tau\nonumber\\
&\leq C+4\max(Nk_1,c)\int_0^t\tau\left[\int_{\mathbb{R}^N}|\nabla h(|u|^2)|^2dx+\frac{1}{2}\int_{\mathbb{R}^N}(|W|*|u|^2)|u|^2dx\right]d\tau.\label{914w2}
\end{align}
Let
$$
A_3(t):=4\int_0^t\tau\left[\int_{\mathbb{R}^N}|\nabla h(|u|^2)|^2dx+\frac{1}{2}\int_{\mathbb{R}^N}(|W|*|u|^2)|u|^2dx\right]d\tau.
$$
(\ref{914w2}) implies
$$
A'_3(t)\leq \frac{C_0}{t}+\frac{\max(Nk_1,c)}{t}A_3(t).
$$
Using Gronwall's inequality, we obtain
$$
A_3(t)\leq t^{\max(Nk_1,c)}[A_3(1)+C-\frac{C}{t^{\max(Nk_1,c)}}]\leq C't^{\max(Nk_1,c)}.
$$
Consequently, we have
$$
\int_{\mathbb{R}^N}|\nabla h(|u|^2)|^2dx+\int_{\mathbb{R}^N}(|W|*|u|^2)|u|^2dx\leq \frac{C}{t^{2-\max(Nk_1,c)}}.
$$

In all cases above, we have
$$
\lim_{t\rightarrow +\infty}\int_{\mathbb{R}^N}|\nabla h(|u|^2)|^2dx+\int_{\mathbb{R}^N}(|W|*|u|^2)|u|^2dx=0.
$$
By the conservation of energy, we get
$$
\lim_{t\rightarrow +\infty}\left(\frac{1}{2}\int_{\mathbb{R}^N}|\nabla u|^2dx+\frac{1}{2}\int_{\mathbb{R}^N}|\nabla h(|u|^2)|^2dx+\frac{1}{4}\int_{\mathbb{R}^N}(|W|*|u|^2)|u|^2dx\right)=E(u_0),
$$
which means that
$$
\lim_{t\rightarrow +\infty}\int_{\mathbb{R}^N}|\nabla u|^2dx=2E(u_0).
$$
By the conservation of mass, we obtain
$$
\lim_{t\rightarrow +\infty}\|u(\cdot,t)\|_{H^1}^2=\lim_{t\rightarrow +\infty}\left(\int_{\mathbb{R}^N}|u|^2dx+\int_{\mathbb{R}^N}|\nabla u|^2dx\right)=M(u_0)+2E(u_0).
$$
(\ref{9181}) is proved.

2. Assume that $u$ is the blowup solution of (\ref{1}),  $u_0\in X$ and $xu_0\in L^2(\mathbb{R}^N)$, $W(x)\geq 0$ and $2W+(x\cdot \nabla W)\leq 0$, $2h''(s)h'(s)s+(h'(s))^2\leq 0$. Using (\ref{893}), we have
\begin{align}
&\quad 2(T-t)^2\int_{\mathbb{R}^N}(W*|u|^2)|u|^2dx\nonumber\\
&=\int_{\mathbb{R}^N}|(x+2i(T-t)\nabla)u|^2dx+4(T-t)^2\int_{\mathbb{R}^N}|\nabla h(|u|^2)|^2dx\nonumber\\
&\quad+4\int_0^t(T-\tau)\theta(\tau)d\tau-8T^2E(u_0)-\int_{\mathbb{R}^N}|xu_0|^2dx\nonumber\\
&\quad-4T\Im\int_{\mathbb{R}^N}\bar{u}_0(x\cdot \nabla u_0)dx.\label{894}
\end{align}
If $$-8T^2E(u_0)-\int_{\mathbb{R}^N}|xu_0|^2dx-4T\Im\int_{\mathbb{R}^N}\bar{u}_0(x\cdot \nabla u_0)dx>0,$$ then (\ref{894}) implies that
\begin{align*}
\int_{\mathbb{R}^N}(W*|u|^2)|u|^2dx\geq \frac{C}{(T-t)^2}.
\end{align*}

Using energy conservation law $E(u)=E(u_0)$, we get
$$
\frac{1}{2}\int_{\mathbb{R}^N}|\nabla u|^2dx+\frac{1}{2}\int_{\mathbb{R}^N}|\nabla h(|u|^2)|^2dx=\frac{1}{4}\int_{\mathbb{R}^N}(W*|u|^2)|u|^2dx+E(u_0)\geq
\frac{C}{(T-t)^2}+E(u_0).
$$
As $t$ close to $T$ enough, we have
$$
\frac{C}{(T-t)^2}+E(u_0)\geq \frac{C'}{(T-t)^2}.
$$
for some constant $0<C'<C$. Hence
$$
\int_{\mathbb{R}^N}|\nabla u|^2dx+\int_{\mathbb{R}^N}|\nabla h(|u|^2)|^2dx\geq \frac{2C'}{(T-t)^2},
$$
(\ref{895}) holds.\hfill $\Box$

{\bf Remark 6.1.} 1. Let $h(s)=s^{\alpha}$ and
\begin{equation}
\label{561} W(x)=m(r)=\left\{
\begin{array}{lll}
-\frac{1}{|x|^m},\quad 0<r\leq 1\\
g(r),\quad 1\leq r\leq 2\\
-\frac{1}{|x|^M},\quad r\geq 2,
\end{array}\right.
\end{equation}
$M>\max(N,m)$, and $g(r)<0$ satisfies $mg(r)+rg'(r)\geq 0\geq Mg(r)+rg'(r)$ if $1\leq r\leq 2$.
Then $2h''(s)h'(s)s+(h'(s))^2\geq 0$ if $\alpha\geq \frac{1}{2}$, while $-(1-2\alpha)(h'(s))^2\leq 2h''(s)h'(s)s+(h'(s))^2\leq 0$ if $0<\alpha\leq\frac{1}{2}$. And $2W+x\cdot \nabla W\geq 0$ if $m\geq 2$, while $-(2-m)|W|\leq 2W+x\cdot \nabla W<0$ if $m<2$.

2. Let $h(s)=s^{\alpha}$ and $W(x)$ be the function defined in (\ref{55w1}). Then $2h''(s)h'(s)s+(h'(s))^2\leq 0$ if $0<\alpha\leq \frac{1}{2}$, $2W+x\cdot \nabla W\leq 0$ if $p>2$.

Therefore, we can obtain the corresponding conclusions of Theorem 5.

\end{document}